\documentclass[11pt,letterpaper]{article}
\usepackage[margin=3.4cm]{geometry}
\usepackage[utf8]{inputenc}
\usepackage{amsmath,amsfonts,amssymb,amsthm,mathrsfs,bm,mathtools,accents,bbm}

\usepackage[hidelinks]{hyperref}
\usepackage[all]{xy}
\usepackage{pgf}
\usepackage{float}
\usepackage{comment}
\usepackage{caption}
\usepackage[british]{babel}

\usepackage[toc,page]{appendix}

\usepackage{xcolor}

\theoremstyle{plain}

\theoremstyle{definition}

\newtheoremstyle{break}{3pt}{3pt}{\itshape}{}{\bfseries}{.}{\newline}{}
\theoremstyle{break}

\DeclareFontFamily{U}{stixextrai}{}
\DeclareFontShape{U}{stixextrai}{m}{n}
 { <-> stix-extra1 }{}

\newcommand*{\parallelogram}{%
  \rlap{\rotatebox{-30}{\rule[.05ex]{.4pt}{.77em}}}%
  \kern.04em%
  \rlap{\kern.36em\raisebox{0.649519052835em}{\rule{.6em}{.4pt}}}%
  \rule{.6em}{.4pt}\kern-.04em%
  \rotatebox{-30}{\rule[.05ex]{.4pt}{.77em}}}

\DeclareMathOperator*{\GL}{GL}

\DeclareMathOperator*{\GO}{GO}
\DeclareMathOperator*{\GSp}{GSp}
\DeclareMathOperator*{\Diag}{Diag}

\newcommand*{\HG}[4]{#1\left(\begin{array}{c|}
#2 \\ #3
\end{array}\,\,
#4
\right)}

\usepackage{tikz}
\usetikzlibrary{tikzmark,arrows,matrix,decorations.pathreplacing,calc,decorations.markings}

\pgfkeys{tikz/mymatrixenv/.style={decoration=brace,every left delimiter/.style={xshift=3pt},every right delimiter/.style={xshift=-3pt}}}
\pgfkeys{tikz/mymatrix/.style={matrix of math nodes,left delimiter={(},right delimiter={)},inner sep=2pt,column sep=1em,row sep=0.5em,nodes={inner sep=0pt}}}
\pgfkeys{tikz/mymatrixbrace/.style={decorate,thick}}

\makeatletter
\newcommand*{\bigominus}{\DOTSB\bigominus@\slimits@}

\newcommand{\bigominus@}{\mathop{\mathpalette\bigominus@@\relax}}

\newcommand{\bigominus@@}[2]{%
  \vcenter{\hbox{%
    \sbox\z@{$\m@th#1\bigoplus$}%
    \resizebox{\wd\z@}{!}{$\m@th#1\bm{\ominus}$}%
  }}%
}
\makeatother

\DeclareFontFamily{U}{mathx}{\hyphenchar\font45}
\DeclareFontShape{U}{mathx}{m}{n}{
      <5> <6> <7> <8> <9> <10>
      <10.95> <12> <14.4> <17.28> <20.74> <24.88>
      mathx10
      }{}
\DeclareSymbolFont{mathx}{U}{mathx}{m}{n}
\DeclareFontSubstitution{U}{mathx}{m}{n}
\DeclareMathAccent{\widecheck}{0}{mathx}{"71}
\DeclareMathAccent{\wideparen}{0}{mathx}{"75}

\newcommand\sidecomment[5][0.3,0.1]%
  {\begin{tikzpicture}[remember picture,overlay]
   \draw[-stealth',thick]
     ($({pic cs:#4}|-{pic cs:#2})+(#1)$)
     .. controls +(1,0) and +(1,0) ..
     node[right,align=left]{#5}
     ($({pic cs:#4}|-{pic cs:#3})+(#1)$);
   \end{tikzpicture}%
  }

\newsavebox{\overlongequation}

\title{Bailey-type factorizations for Horn functions.}
\author{Carlo Verschoor}

\AtEndDocument{\bigskip{\footnotesize
  \noindent
  \textsc{Carlo Verschoor, Utrecht University}
  \textit{E-mail address}: \texttt{\href{mailto:carlovmm@gmail.com}{carlovmm@gmail.com}}
}}

\begin{document}

\maketitle

\begin{abstract}
In this paper we are interested in extending Bailey's identity to other classical hypergeometric functions. Bailey's identity states that under a suitable choice of parameters, Appell's $F_4$ decomposes into a product of two ${}_2F_1$'s. We will show how Bailey-type factorizations can be found for Horn's hypergeometric functions $H_1, H_4$ and $H_5$.
\end{abstract}

\section{Introduction}
\label{chap-relations}

In 1933, Bailey published the following identity \cite{Bai33} between classical hypergeometric functions

\begin{align}
    \notag
    &\HG{F_4}{a,\,\,b}{c,\,\,a+b-c+1}{x(1-y),y(1-x)} \\&= \HG{{}_2 F_1}{a,\,\,b}{c}{x}\HG{{}_2 F_1}{a,\,\,b}{a+b-c+1}{y}. \label{eq-bailey-F4}
\end{align}

Here ${}_2F_1$ is Gauss' hypergeometric function defined as
$$\HG{{}_2F_1}{a,\,\,b}{c}{x} = \sum_{n=0}^{\infty} \frac{(a)_n(b)_n}{(c)_nn!}z^n.$$

And $F_4$ is Appell's $F_4$ hypergeometric function defined as
$$\HG{F_4}{a,\,\,b}{c,\,\,c'}{x,y} = \sum_{m=0}\sum_{n=0} \frac{(a)_{m+n}(b)_{m+n}}{(c)_m(c')_{n}m!n!}x^my^n.$$

Other identities similar to \eqref{eq-bailey-F4} for Appell's $F_2$ have been found independently by Beukers \cite{Zud16} and Vidunas \cite{Vid08}. It reads,

{\footnotesize
    $$\HG{F_2}{a+b-\frac{1}{2},\,\,a,\,\,b}{2a,\,\,2b}{\frac{4u(1-u)(1-2v)}{(1-2uv)^2},\frac{4v(1-v)(1-2u)}{(1-2uv)^2}}$$ 
    \begin{equation}=(1-2uv)^{-1+2a+2b}\HG{{}_2 F_1}{a+b-\frac{1}{2},\,\,a}{2a}{4u(1-u)}\HG{{}_2 F_1}{a+b-\frac{1}{2},\,\,b}{2b}{4v(1-v)}.\label{eq-bailey-F2}\end{equation}
}

For similar identities related to Appell's $F_3$ it is noted in \cite{Vid08} that these follow directly from the observation that the two functions

\begin{equation*}
\HG{F_2}{a\,\, b_1,\,\,b_2}{c_1,\,\,c_2}{x,y},
\end{equation*}
\begin{equation*}
x^{-b_1}y^{-b_2}\HG{F_3}{1+b_1-c_1\,\, 1+b_2-c_2,\,\,,b_1,\,\,b_2}{1+b_1+b_2-a}{\frac{1}{x},\frac{1}{y}}
\end{equation*}
satisfy the same differential equations.

The aim of this paper is to generalize \eqref{eq-bailey-F4} 
to other classical two-variable hypergeometric functions.

Let ${}_{2,2}F_{1,1}(x,y)$ be such a function. Then we look for factorizations of the
form
\begin{equation}\label{baileyfactorization}
{}_{2,2}F_{1,1}(\phi(s,t),\psi(s,t))=Q(s,t)^\lambda
\HG{{}_2 F_1}{a,\,\,b}{c}{\rho(s)}\HG{{}_2 F_1}{a',\,\,b'}{c'}{\sigma(t)},
\end{equation}
where $\phi,\psi,Q$ are rational functions in two variables and $\rho,\sigma$ rational
functions in one variable. Analytic continuation of the right hand side in $s,t$
will generate a space of functions of dimension $2\times2=4$. Hence we expect the function
on the left to be a solution of a rank $4$ equation. Which is the case for Appell
$F_2,F_4$.

In this paper we are interested in Bailey-type factorizations for Horn functions of rank $4$, in particular $H_1, H_4$ and $H_5$. These are defined as the hypergeometric series

\begin{align}
    \HG{H_1}{a,\,\,b,\,\,c}{d}{x,y} &= \sum_{m=0}^{ \infty}\sum_{n=0}^{\infty}\frac{(a)_{m-n}(b)_{m+n}(c)_n}{(d)_m m! n!} x^m y^n, \label{eq-H1} \\
    \HG{H_4}{a,\,\,b}{c,\,\,d}{x,y} &= \sum_{m=0}^{ \infty}\sum_{n=0}^{\infty}\frac{(a)_{2m+n}(b)_{n}}{(c)_m(d)_n m! n!} x^m y^n, \label{eq-H4} \\
    \HG{H_5}{a,\,\,b}{c}{x,y} &= \sum_{m=0}^{ \infty}\sum_{n=0}^{\infty}\frac{(a)_{2m+n}(b)_{n-m}}{(c)_n m! n!} x^m y^n. \label{eq-H5}
\end{align}

\section{Motivation}
We can see why hypergeometric functions factor as Bailey-type identities, by looking at monodromy. 
The monodromy group of ${}_2 F_1(x)$ is contained in
$\GL(2)$. Hence the monodromy group of ${}_{2,2}F_{1,1}$ is contained in
$\GL(2) \times \GL(2)$. Over the complex numbers
the group $\GL(2)\times\GL(2)$ maps injectively into the group of orthogonal similitudes $\GO(4)$ by the map $$(M_1,M_2) \mapsto \left(\begin{pmatrix}a&b\\c&d\end{pmatrix} \mapsto M_1^{\intercal}\begin{pmatrix}a&b\\c&d\end{pmatrix}M_2\right).$$ 

The $\GO(4)$ semi-invariant form is given by the determinant $ad-bc$. 
Thus to find Bailey-like factorizations we should look for a rank $4$ system whose monodromy
is contained in $\GO(4)$. To find candidates for such hypergeometric systems, we
write down a basis of monodromy matrices $M_1,\ldots,M_k$ as constructed by Beukers' method for finding monodromy in \cite{Beu13}. 
Let $Q$ be a candidate non-singular $4\times 4$ matrix and test whether
$M_i^TQM_i=\lambda_iQ$ for some scalar $\lambda_i$. The requirement that we should find
a non-trivial $Q$ gives us a large set
of restrictions on the hypergeometric parameters $\alpha_1,\ldots,\alpha_r$, or more
correctly, $e^{2\pi i\alpha_1},\ldots,e^{2\pi i\alpha_r}$. If the group generated by
$M_1,\ldots,M_k$ acts irreducibly, the resulting matrix $Q$ will be either symmetric
or anti-symmetric. This follows from the fact that $Q$ is uniquely determined up to
a scalar and both $Q$ and $Q^T$ are a solution. In the first case $\langle M_1,\ldots,M_k\rangle$ will
be contained in $\GO(4)$, in the second case in $\GSp(4)$, the group of symplectic
similitudes.

First we give a short summary of how the algorithm works. The implementation follows straightforwardly from the algorithm, so we will omit that part.

The monodromy matrices given to us by Beukers' method in \cite{Beu13} have the form 
$$M_{\bm{\rho},j} := X_{\bm{\rho}} \chi_{\bm{\rho},j} X_{\bm{\rho}}^{-1}$$ where $\chi_{\bm{\rho},j}$ is a diagonal matrix and both $X_{\bm{\rho}}$ and $\chi_{\bm{\rho},j}$ only depend on the parameters $e^{2\pi i\alpha_1},\ldots,e^{2\pi i\alpha_r}$. Here the numbers $\alpha_1,\ldots,\alpha_r$ depends linearly on the classical parameters of the rank $4$ hypergeometric function we want to test.

Now we want to find a (skew)-symmetric matrix $Q$ and a scalar 
$\lambda_{\bm{\rho},j} \in \mathbb{C}^*$ such that 
$$M_{\bm{\rho},j}^{\intercal}Q M_{\bm{\rho},j} = \lambda_{\bm{\rho},j}Q.$$
This can then be written as
$$(X_{\bm{\rho}}^{\intercal})^{-1} \chi_{\bm{\rho},j} X_{\bm{\rho}}^{\intercal} Q X_{\bm{\rho}} \chi_{\bm{\rho},j} X_{\bm{\rho}}^{-1} = \lambda_{\bm{\rho},j}Q.$$
And then it can be written as
$$\chi_{\bm{\rho},j} X_{\bm{\rho}}^{\intercal} Q X_{\bm{\rho}} \chi_{\bm{\rho},j} = \lambda_{\bm{\rho},j} X_{\bm{\rho}}^{\intercal}Q X_{\bm{\rho}}.$$
So let $Q_{\bm{\rho}} = X_{\bm{\rho}}^{\intercal}Q X_{\bm{\rho}}$, then we have 
$$\chi_{\bm{\rho},j} Q_{\bm{\rho}} \chi_{\bm{\rho}_j} = \lambda_{\bm{\rho},j} Q_{\bm{\rho}}.$$
So let $\chi_{\bm{\rho},j} = \Diag(\chi_1,\ldots,\chi_D)$ then for each row-column entry $(r,c)$ where $Q_{\bm{\rho}}$ is non-zero, we need to have $\chi_r \chi_c = \lambda_{\bm{\rho},j}$. 

These conditions result in restrictions on the parameters $\alpha_1,\ldots,\alpha_r$, or more correctly $e^{2\pi i \alpha_1},\ldots,e^{2\pi i \alpha_r}$. 

This turns out to be a reasonably fast approach to
compute which specializations give us a monodromy group contained in $\GO(4)$ or $\GSp(4)$ for the classical cases. We need to be careful though, because the restrictions hold for $e^{2\pi i \alpha_1},\ldots,e^{2\pi i \alpha_r}$ we find $\bm{\alpha}$-vectors upto shifts in $\mathbb{Z}^r$. Table \ref{tab:sometab} gives the results of our implementation of the 
algorithm described. We feed it an $A$-hypergeometric system and it computes for which parameter vectors $\bm{\alpha}$ the system is either in 
$\mbox{GO}(4)$ or $\mbox{GSp}(4)$. The specialization given in Table \ref{tab:sometab} correspond to the classical parameters and so not to the $A$-hypergeometric $\bm\alpha$-vector. Additionally note that $\bm{\alpha}$ may not be totally non-resonant, this is given in the third column. Here TNR means Totally Non-Resonant and NR means Non-Resonant. Note that this may give us an inconsistency, because the
computation of the monodromy matrices was made under the assumption of TNR. However,
NR occurs only in the symplectic case, which will not be considered any further. The fourth column tells us whether the monodromy for the corresponding specialization is contained in $\GO(4)$ or $\GSp(4)$.

Among the classical two-variable hypergeometric equations of rank $4$ the cases
$F_3,H_2,H_7$ are missing from the table. However, $F_3$ and $H_2$ are related to $F_2$
(they have the same $A$-polytope) and $H_7$ is related to $H_4$.

    

\begin{minipage}{\linewidth}
    \bigskip
    \centering
    \captionof{table}{} \label{tab:sometab}
    \begin{tabular}{cc|cc}
    System & Specialization $\pmod{\mathbb{Z}^r}$& Res. & Mon. \\ \hline
    $F_2$ & $(q_0,q_1,q_0-q_1+\frac{1}{2},2q_1,2q_0-2q_1)$ & TNR & $\mbox{GO}(4)$ \\
    & $(q_0,q_1,q_0-q_1,2q_1,2q_0-2q_1)$ & NR & $\mbox{GSp}(4)$ \\&&& \\
    $F_4$ & $(q_0, q_1, q_2, q_0 + q_1 - q_2)$ & TNR & $\mbox{GO}(4)$ \\&&& \\
    
    $H_1$ & $(q+1,q+\frac{1}{2},\frac{1}{2},2q)$ & TNR & $\mbox{GO}(4)$ \\
    & $(q,q,-\frac{1}{2},2q)$ & NR & $\mbox{GSp}(4)$ \\&&& \\
    $H_4$ & $(q_0, q_1, q_0 - q_1, 2 q_1)$ & TNR & $\mbox{GO}(4)$ \\
    & $(q_0 - \frac{1}{2}, q_1, q_0 - q_1,2 q_1)$ & TNR & $\mbox{GSp}(4)$ \\&&& \\
    $H_5$ & $(q+\frac{1}{2},q,2q)$ & TNR & $\mbox{GO}(4)$ \\
    & $(q,q,2q)$ & NR & $\mbox{GSp}(4)$ \\
    \end{tabular}
    \bigskip
\end{minipage}

Looking at $F_4$ and Bailey's identity \eqref{eq-bailey-F4} we see that the parameter vectors 
$(q_0, q_1, q_2, q_0 + q_1 - q_2)$ and its classical parameter vector $(a,b,c,a+b-c+1)$ are the same modulo $\mathbb{Z}^4$ if we take $q_0 = a, q_1 = b, q_2 = c$. The same can be said for Bailey's identity for $F_2$ in \eqref{eq-bailey-F2}. Looking at this table one could now wonder whether $H_1$, $H_4$ and $H_5$ also admit a Bailey type decomposition.

To find a Bailey type factorization we follow an approach from \cite{Vid08} who finds 
such identities for $F_2$ and $F_4$.
Consider \eqref{baileyfactorization} and write $F(x,y)={}_{2,2}F_{1,1}(x,y)$.
Fix $s$ and note that $F(\phi(s,t),\psi(s,t))$, as a function of $t$, satisfies
an ordinary second order differential equation with rational function
coefficients, because the right hand side of \eqref{baileyfactorization} does.
The problem we like to solve is to find (rational) functions $x(t),y(t)$ such
that $f(t):=F(x(t),y(t))$ satisfies an ordinary second order differential equation.
For general $x(t),y(t)$ such a function would satisfy a fourth order equation,
so the second order restriction does give us restrictions on $x(t),y(t)$.

Suppose $f(t)$ is the solution to a second order differential system of the form
$$\frac{d^2f}{dt^2} + c_1 \frac{df}{dt} + c_2f = 0.$$

To find the relation with $F$ and its differential equations we apply
the chain rule and product rule multiple times to $f(t) = F(x(t),y(t))$.
\begin{align}
    \frac{d f}{d t} &= \dot{x}F_x + \dot{y} F_y \label{eq:Ft} \\ 
    \frac{d^2 f}{d t^2} &= \dot{x}^2 F_{xx} + 2\dot{x}\dot{y} F_{xy} + \dot{y}^2F_{yy} + \ddot{x}F_{x} + \ddot{y}F_y \label{eq:Ftt}
\end{align}
Here we denote $$\dot{x} = \frac{dx}{dt},\,\,\,\dot{y} = \frac{dy}{dt}.$$

\section{Horn's $H_4$}
\label{chap-H4}
The first system we want to investigate is Horn's $H_4$ hypergeometric function \eqref{eq-H4}. This one is interesting because we can see from Table \ref{tab:sometab} that it's corresponding system has a monodromy group in $\GO(4)$ and the corresponding specialization of the parameters is two dimensional. 
A system of partial differential equations for $H_4$ can be found at \cite[p.817]{Deb02}. It is
\begin{align}&x(1 - 4x)F_{xx} -4xy F_{xy} - y^2F_{yy} \nonumber \\
&\,\,\,+ (c - (4a + 6)x)F_x - 2(a + 1)yF_y - a(a + 1)F = 0, \label{eq:H4_1} \\
&y (1 - y )F_{yy} - 2xy F_{xy} \nonumber \\
&\,\,\,+ (d - (a + b + 1)y )F_y - 2bxF_x - abF = 0. \label{eq:H4_2}
\end{align}


Now we would like to parameterize $x$ and $y$ with variable $t$ such that we get an equation of the form 
$$\frac{d^2 f}{d t^2} + c_1 \frac{d f}{d t} + c_2 f = 0.$$
To achieve this we take \eqref{eq:Ftt} and eliminate
partial derivatives on the righthand side terms of \eqref{eq:Ftt}.
First we can eliminate $F_{xx}$ and
$F_{yy}$ 
using equations \eqref{eq:H4_1} and \eqref{eq:H4_2}. Then
we set the coefficient of 
$F_{xy}$ to be zero. This coefficient is equal to
\begin{equation}\frac{-2y(y - 2)\dot{x}^2 + 2(y-1)(4x - 1)\dot{x}\dot{y} - 2x(4x - 1) \dot{y}^2}{(y-1)(4x-1)}.\label{eq:H4_coef}\end{equation}

If we specialize this by
$$(x(u,v),y(u,v)) = \left(\frac{1 - \left(\frac{v - \frac{1}{v}}{2}\frac{u - \frac{1}{u}}{2}\right)^2}{4},1 + \frac{v + \frac{1}{v}}{2}\right),$$
we can factor \eqref{eq:H4_coef} as
$$\frac{(v^2\dot{u} + 2u\dot{v} - \dot{u})(v^2\dot{u} - 2u\dot{v} - \dot{u})(u^2 + 1)^2(v^2 - 1)^2}{64(v^2 + 1)u^4v^3}.$$
This is equal to zero when $v= \pm 1$ or $u = \pm i$, but we skip those cases and focus on one of the two differential equations that emerge (the other one gives a similar result). Thus consider the following differential equation:
$$(v^2 - 1) \dot{u} = -2u \dot{v}.$$
Which can be simplified to
$$\frac{-1}{2u} du = \frac{1}{v^2 - 1} dv.$$
After integration we obtain
$$-\frac{1}{2}\log(u) + C = \frac{1}{2}(\log(1-v) - \log(1+v)).$$
And thus we get the solutions
$$u =C' \frac{1+v}{1-v},\mbox{ $C'$ constant}.$$
Now set $v=t$, $u = C'\frac{1 + t}{1 - t}$ and let $C' = \frac{1 + s}{1 - s}$ with $s$ another constant. 
This gives us a parameterization
$$x(s,t) = -\frac{(st^2 + s + 2t)(2st + t^2 + 1)s}{4(s + 1)^2(s - 1)^2t^2},\,\,\,\, y(s,t) = 1 + \frac{t + \frac{1}{t}}{2},$$
which ensures that \eqref{eq:H4_coef} is $0$.

Now we have eliminated $F_{xx},F_{xy},F_{yy}$ from \eqref{eq:Ftt}. 
Equation \eqref{eq:Ftt} has now acquired the form
\begin{equation}\frac{d^2f}{dt^2}=c_3F_x(x(t),y(t))+c_4F_y(x(t),y(t))+c_5f.\label{eq-c3c4c5-form}\end{equation}
The coefficients $c_3,c_4,c_5$ are quite cumbersome to write down and can be found in Appendix \ref{sec-app-h4}.

Then, what we would like to see is that
the vector $(c_3,c_4)$ is a multiple, say $r_2(t)$, of $(\dot{x},\dot{y})$.
In such a case we have $$c_3F_x+c_4F_y=r_2(t)(\dot{x}F_x+\dot{y}F_y)=r_2\frac{df}{dt}$$
in virtue of equation \eqref{eq:Ft}.
It is a small miracle that this indeed happens for the
parameter choices $a = q_0, b = q_1, c = 1 + q_0 - q_1, d = 2q_1$. This choice has been motivated by Table \ref{tab:sometab}.
This means that for fixed $s$ the function
$$\HG{H_4}{q_0,\,\, q_1}{1 + q_0 - q_1,\,\, 2q_1}{-\frac{{\left(s t^{2} + s + 2 \, t\right)} {\left(2 \, s t + t^{2} + 1\right)} s}{4 \,{\left(s^2 - 1\right)}^{2} t^{2}},\frac{{\left(t + 1\right)}^{2}}{2 \, t}}.$$
satisfies an ordinary Fuchsian equation of the form
\begin{equation}
    \frac{d^2F}{dt^2} + r_2 \frac{dF}{dt} + r_3 F \label{eq:Fuchs2}.
\end{equation}
Here $r_2$ and $r_3$ are some rational coefficients in $q_0,q_1,s$ and $t$ which can be found in Appendix \ref{sec-app-h4}.
It is a Fuchsian equation which we like to identify with a transform of the Gaussian
hypergeometric equation. To that end we investigate its local exponents.

The following singularities with corresponding local exponents were found:
\begin{center}
\begin{tabular}{c|cc}
Singularity & Exponent 1 & Exponent 2 \\ \hline
$t = 1$ & $0$ & $2$ \\
$t = -1$ & $0$ & $2 - 4q_1$ \\
$t = 0$ & $q_0$ & $q_0 + 1$ \\
$t = \infty$ & $q_0$ & $q_0 + 1$ \\
Roots of $t^2 + 2t/s + 1 = 0$ & $0$ & $q_1 - q_0$ \\
Roots of $t^2 + 2ts +  1 = 0$ & $0$ & $q_1 - q_0$
\end{tabular}
\end{center}
At $t=0$ and $t=\infty$ we notice there is a difference of $1$ in the local exponents, hence they may be apparent singularities. As the system for ${}_2 F_{1}$ only has local exponents at 
$0,1$ and $\infty$, we might want to map the differential equations for ${}_2 F_{1}(z)$ under the covering $$z \mapsto \frac{(s+1)^2}{(s - 1)^2}\frac{(t+1)^4}{(t - 1)^4}.$$
This is because $t = 1$ corresponds to $z = \infty$ with multiplicity $4$; And $t = -1$ corresponds to $z = 0$ with multiplicity $4$; If $t$ is a root of $t^2 + 2t/s + 1 = 0$ or $t^2 + 2ts +  1 = 0$ this corresponds to $z = 1$ each with multiplicity $1$. And $t=0$ and $t = \infty$ correspond to regular points. 

Recall that the Riemann scheme of $\HG{{}_2 F_1}{a,\,\,b}{c}{z}$ is equal to 
\begin{center}
\begin{tabular}{c|cc}
Singularity & Exponent 1 & Exponent 2 \\ \hline
$z = 0$ & $0$ & $1-c$ \\
$z = 1$ & $0$ & $c-a-b$ \\
$z = \infty$ & $a$ & $b$ \\
\end{tabular}
\end{center}

Hence the Riemann scheme of $\HG{{}_2 F_1}{a,\,\,b}{c}{\frac{(s+1)^2}{(s - 1)^2}\frac{(t+1)^4}{(t - 1)^4}}$ becomes

\begin{center}
\begin{tabular}{c|cc}
Singularity & Exponent 1 & Exponent 2 \\ \hline
$t = 1$ & $4a$ & $4b$ \\
$t = -1$ & $0$ & $4 - 4c$ \\
$t = 0$ & $0$ & $1$ \\
$t = \infty$ & $0$ & $1$ \\
Roots of $t^2 + 2t/s + 1 = 0$ & $0$ & $c-a-b$ \\
Roots of $t^2 + 2ts +  1 = 0$ & $0$ & $c-a-b$
\end{tabular}
\end{center}
Upto translation these exponents should be equal to the local exponents we found for $H_4$, hence the differences between local exponents form a linear set of equations. $4b - 4a = 2$, $4 - 4c = 2 - 4q_1$ and $c - a - b = q_1 - q_0$. This can be solved for $a = \frac{1}{2}q_0 + \frac{1}{2}$, $b = \frac{1}{2}q_0$ and $c = q_1 + \frac{1}{2}$. Indeed we can check that under these transformations from the ${}_2 F_1$ system we obtain equation \eqref{eq:Fuchs2}.

Now we need to translate the local exponents. For this we need to multiply by an additional factor of $(1-t)^{-2q_0}t^{q_0}$. This makes the local exponents at $t = 0$ and $t = \infty$ both $q_0$ and $q_0 + 1$ and the local exponents at $t = 1$ now become $0$ and $2$.

We conclude that
$$\HG{H_4}{q_0,\,\, q_1}{1 + q_0 - q_1,\,\, 2q_1}{-\frac{{\left(s t^{2} + s + 2 \, t\right)} {\left(2 \, s t + t^{2} + 1\right)} s}{4 \, {\left(s^2 - 1\right)}^{2} t^{2}},\frac{{\left(t + 1\right)}^{2}}{2 \, t}}$$
satisfies the same second order ordinary differential equation in $t$ as
$$(1-t)^{-2q_0} t^{q_0} \HG{{}_2 F_1}{\frac{1}{2}q_0 + \frac{1}{2},\,\, \frac{1}{2}q_0}{q_1 + \frac{1}{2}}{\frac{(s+1)^2(t+1)^4}{(s-1)^2(t-1)^4}}.$$

Now we want to transform this such that it becomes symmetric in the arguments $s$ and $t$ in the sense that swapping $s$ and $t$ does not change the system, this can be done with the following transformation:
$$s \to -\frac{s^2 + 1}{s^2 - 1},\,\,\, t \to
 \frac{st + 1}{st - 1}.$$
In this way we obtain
\begin{equation}\HG{H_4}{q_0,\,\, q_1}{1 + q_0 - q_1,\,\, 2q_1}{\frac{{\left(s^4 - 1\right)}
{\left(t^4 - 1\right)}}{4 \,{\left(s^2t^2 - 1\right)}^{2}},
\frac{2s^2t^2}{s^2t^2 - 1}}. \label{eq:H4_dif_left}\end{equation}
And as $s$ is a constant it should satisfy the same differential equation in $t$ as
{\small 
\begin{equation}\left(\frac{1}{1-st}\right)^{-2q_0} \left(\frac{st+1}{st - 1}\right)^{q_0}
\HG{{}_2 F_1}{\frac{1}{2}q_0 + \frac{1}{2},\,\, \frac{1}{2}q_0}{q_1 + \frac{1}{2}}{s^4}\HG{{}_2 F_1}{\frac{1}{2}q_0 + \frac{1}{2},\,\, \frac{1}{2}q_0}{q_1 + \frac{1}{2}}{t^4}.\label{eq:H4_dif_right}\end{equation}}

By the symmetry in $s,t$ we can also say that \eqref{eq:H4_dif_left} and 
\eqref{eq:H4_dif_right} satisfy the same second order ordinary differential equation in $s$.

To establish a Bailey-like identity we consider \eqref{eq:H4_dif_left} near the
point $(s,t)=(1,0)$. It is holomorphic there, so \eqref{eq:H4_dif_left} is $(1-s^2t^2)^{q_0}$ times a hypergeometric series holomorphic near $s=1$ times a hypergeometric series 
holomorphic near $t=0$.

After setting the constant terms both equal to $1$ we conclude that
\begin{equation*}\HG{H_4}{q_0,\,\, q_1}{1 + q_0 - q_1,\,\, 2q_1}{\frac{{\left(s^4 - 1\right)}{\left(t^4 - 1\right)}}{4 \,{\left(s^2t^2 - 1\right)}^{2}},
\frac{2s^2t^2}{s^2t^2 - 1}}
\end{equation*}
\begin{equation*}
=\left(1 - s^2t^2\right)^{q_0}
\HG{{}_2 F_1}{\frac{1}{2}q_0 + \frac{1}{2},\,\, \frac{1}{2}q_0}{q_0 - q_1 + 1}{1 - s^{4}}\HG{{}_2 F_1}{\frac{1}{2}q_0 + \frac{1}{2},\,\, \frac{1}{2}q_0}{q_1 + \frac{1}{2}}{t^{4}}.
\end{equation*}
Notice that only squared variables are used, so substitute $s \to \sqrt{s}$ and $t \to \sqrt{t}$ to obtain:

\begin{equation*}\HG{H_4}{q_0,\,\, q_1}{1 + q_0 - q_1,\,\, 2q_1}{\frac{{\left(s^2 - 1\right)}{\left(t^2 - 1\right)}}{4 \,{\left(st - 1\right)}^{2}},
\frac{2st}{st - 1}}
\end{equation*}
\begin{equation}
=\left(1 - st\right)^{q_0}
\HG{{}_2 F_1}{\frac{1}{2}q_0 + \frac{1}{2},\,\, \frac{1}{2}q_0}{q_0 - q_1 + 1}{1 - s^{2}}\HG{{}_2 F_1}{\frac{1}{2}q_0 + \frac{1}{2},\,\, \frac{1}{2}q_0}{q_1 + \frac{1}{2}}{t^{2}}. \label{eq:H4_rel}
\end{equation}

After a search in the literature it turned out that there is another way
to get the same identity. In \cite[eq. 7.6, p. 382]{Erd48} we find the identity
$$
\HG{H_4}{\alpha,\beta}{\gamma,2\beta}{x,y}=(1-y/2)^{-\alpha}
\HG{F_4}{\frac12\alpha,\frac12\alpha+\frac12}{\gamma,\beta+\frac12}
{\frac{16x}{(2-y)^2},\frac{y^2}{(2-y)^2}}.
$$
Using this identity we obtain
\begin{align*}
&\HG{H_4}{q_0,q_1}{1+q_0-q_1,2q_1}{\frac{(s^2-1)(t^2-1)}{4(st-1)^2},\frac{2st}{st-1}}\\
&=(1-st)^{q_0}\HG{F_4}{\frac12 q_0,\frac12 q_0+\frac12}{1+q_0-q_1,q_1+\frac12}
{(1-s^2)(1-t^2),s^2t^2}.
\end{align*}
We can now apply Bailey's factorization \eqref{eq-bailey-F4} to get the right hand
side of \eqref{eq:H4_rel}.

\section{Horn's $H_1$}
\label{chap-H1}

A system of partial differential equations for Horn's $H_1$ function \eqref{eq-H1} can be found at \cite[p.817]{Deb02}. It is given by
\begin{align}&x(1 - x)F_{xx}+ y^2F_{yy} \nonumber \\
&\,\,\,+ (d - (a+b+1)x)F_x - (a-b-1)yF_y - abF = 0, \label{eq:H1_1} \\
&y (1 + y )F_{yy} - x(1-y) F_{xy} \nonumber \\
&\,\,\,+ (a-1 - (b + c + 1)y )F_y - cxF_x - bcF = 0. \label{eq:H1_2}
\end{align}
We can use \eqref{eq:H1_1} and \eqref{eq:H1_2} to eliminate $F_{xx}$ and $F_{xy}$ from \eqref{eq:Ftt}. The following coefficient for $F_{yy}$ remains:

\begin{equation}
\frac{(y^3 - y^2)\dot{x}^2 + (-2xy^2 - 2xy + 2y^2 + 2y)\dot{x}\dot{y} + (x^2y - x^2 - xy + x)\dot{y}^2}{(x-1)x(y-1)}. \label{eq:H1_coef}
\end{equation}

Using the specialization
$$x(u,v) = 1 - \frac{(v - 1/v)^2(u-1/u)^2}{16},\,\,\,y(u,v) = v^2,$$
the coefficient \eqref{eq:H1_coef} factors as follows
$$\frac{4(v^2\dot{u} + 2u\dot{v} - \dot{u})(v^2\dot{u} - 2u\dot{v} - \dot{u})(u^2 + 1)^2v^4}{((uv + u + v - 1)(uv + u - v + 1)(uv - u + v + 1)(uv - u - v - 1)u^2)}.$$
Just like in the $H_4$ case we can make this zero by taking the parametrization
$$u = C \frac{1 + t}{1-t},\quad v=t.$$
And again we will pick the integration constant $C$ to be
$$C = \frac{1 + s}{1-s}.$$

Now in equation \eqref{eq:Ftt} we eliminated the coefficients $F_{xx}, F_{xy}$ and $F_{yy}$, and we are left to eliminate $F_x$ and $F_y$ using \eqref{eq:Ft}. Call the coefficients for $F_x$ and $F_y$ respectively $c_3$ and $c_4$ and the constant coefficient $c_5$. These coefficients can be found in Appendix \ref{sec-app-h1}. Then what we would like to see is that $c_4 / c_3 = \dot{y} /\dot{x}$. Again miraculously this
happens when $a = q_0 - \frac{1}{2}, b = q_0, c = \frac{1}{2}, d = 2q_0$, which matches the entry in Table \ref{tab:sometab} modulo $\mathbb{Z}^4$. This means that for fixed $s$ the function
$$\HG{H_1}{q_0 - \frac{1}{2},\,\, q_0, \,\, \frac{1}{2}}{2q_0}{-\frac{{\left(s t^{2} + s + 2 \, t\right)} {\left(2 \, s t + t^{2} + 1\right)} s}{{\left(s^2 - 1\right)}^{2} t^{2}},t^2}$$
satisfies an ordinary Fuchsian equation of the form
$$\frac{d^2F}{dt^2} + r_2 \frac{dF}{dt} + r_3 F.$$
The coefficients $r_2$ and $r_3$ can be found at Appendix \ref{sec-app-h1}. Let us analyze the local exponents of this differential equation.

\begin{center}
\begin{tabular}{c|cc}
Singularity & Exponent 1 & Exponent 2 \\ \hline
$t = 0$ & $2q_0 - 1$ & $2q_0$ \\
$t = \infty$ & $2q_0$ & $2q_0 + 1$ \\
Roots of $t^2 + 2t/s + 1 = 0$ & $0$ & $1 - 2q_0$ \\
Roots of $t^2 + 2ts +  1 = 0$ & $0$ & $1 - 2q_0$
\end{tabular}
\end{center}
We see that in this case $t=0$ and $t=\infty$ may be apparent singularities, but $t=1$ and $t=-1$ are regular. Assume just like in the $H_4$ case that we can go from ${}_2 F_1$ to this by the covering $z \mapsto \frac{(s+1)^2(t+1)^4}{(s-1)^2(t-1)^4}$ now we still want $t=1$ ($z=\infty$ with multiplicity $4$) and $t=-1$ ($z=0$ with multiplicity $4$) to be regular. 
Hence $1 - c = \frac{1}{4}$, so $c=\frac{3}{4}$. And $c-a-b = 1 - 2q_0$ and $b - a = \frac{1}{4}$ so $a = q_0 - \frac{1}{4}$ and $b = q_0$.
Now we need to translate the local exponents such that $t=0$ and $t = \infty$ become apparent singularities and $t=1$ becomes regular, for this we need to multiply by the function by $(t-1)^{1-4q_0}t^{2q_0 - 1}$

We obtain that
$$\HG{H_1}{q_0 - \frac{1}{2},\,\, q_0, \,\, \frac{1}{2}}{2q_0}{-\frac{{\left(s t^{2} + s + 2 \, t\right)} {\left(2 \, s t + t^{2} + 1\right)} s}{{\left(s^2 - 1\right)}^{2} t^{2}},t^2}$$

satisfies the same second order differential equation in $t$ as
$$(t-1)^{1-4q_0}t^{2q_0 - 1}\HG{{}_2 F_1}{q_0 - \frac{1}{4},\,\, q_0}{\frac{3}{4}}{\frac{(s+1)^2(t+1)^4}{(s-1)^2(t-1)^4}}.$$

Now we want to make the arguments symmetric again. We can do this by the substitution
$$s \to -\frac{s^2 + 1}{s^2 - 1},\,\,\, t \to
\frac{st + 1}{st - 1}.$$

As a consequence
$$\HG{H_1}{q_0 - \frac{1}{2},\,\, q_0, \,\, \frac{1}{2}}{2q_0}{
\frac{(s^4 - 1) (t^4 - 1)}{(s^2t^2 - 1)^2},
\left(\frac{st + 1}{st - 1}\right)^2
}$$
satisfies the same differential equation in $t$ as
$$(st-1)^{4q_0 - 1}\left(\frac{st + 1}{st - 1}\right)^{2q_0 - 1} 
\HG{{}_2 F_1}{q_0 - \frac{1}{4},\,\, q_0}{\frac{3}{4}}{s^4}
\HG{{}_2 F_1}{q_0 - \frac{1}{4},\,\, q_0}{\frac{3}{4}}{t^4}.$$
Hence by symmetry this is true both in $s$ and $t$.
To find a Bailey-like identity, which is an identity of two variable
power series, we must look at the function
$$\phi: (s,t)\mapsto \left(\frac{(s^4 - 1) (t^4 - 1)}{(s^2t^2 - 1)^2},
\left(\frac{st + 1}{st - 1}\right)^2\right).$$
Then we wish to choose a point $(s_0,t_0)$ with $s_0^4,t_0^4\in\{0,1\}$ such
that $\phi$ maps an open neighbourhood of $(s_0,t_0)$ to a neighbourhood of $(0,0)$.
Unfortunately this is impossible. Such points must satisfy $s_0^4=t_0^4=1$ and $s_0t_0=-1$.
For example, $s_0=-1,t_0=1$. But $\phi$ is not continuous in $(-1,1)$.
In order to get a meaningful identity we could restrict $\phi$ to a neighbourhood $U$
of $(-1,1)$ of the form $s=-1+2uv,t=1+2v$ with $u,v$ small. One verifies that
$$\phi(-1+2uv,1+2v)=(-4u+\mbox{h.o.t.}\,,\,v^2+\mbox{h.o.t.}),$$
where 'h.o.t.' means 'higher order terms'. So $\phi$ is well-defined on $U$ and
its image is in a neighbourhood of $(0,0)$.

Therefore
\begin{equation}
\HG{H_1}{q_0 - \frac{1}{2},\,\, q_0, \,\, \frac{1}{2}}{2q_0}{\phi(-1+2uv,1+2v)}{}
\label{eq:H1leftside}
\end{equation}
is holomorphic near the point $(u,v)=(0,0)$. The only function in the space spanned
by the products of Gauss hypergeometric functions can be
\begin{align*}
&(st-1)^{2q_0}(st+1)^{2q_0 - 1} (1-t^4)^{1 - 2q_0}\\
&\times\HG{{}_2 F_1}{1-q_0,\,\, \frac{3}{4}-q_0}{2 - 2q_0}{1-t^4}
\HG{{}_2 F_1}{q_0 - \frac{1}{4},\,\, q_0}{2q_0}{1-s^4}.
\end{align*}
After the substitution $s\to-1+2uv,t\to1+2v$ this becomes, after normalization,
\begin{align}
&(1+v-uv-2uv^2)^{2q_0}\left(\frac{1-u-2uv}{(1+v)(1+2v+2v^2)}\right)^{2q_0 - 1}
\nonumber \\
&\quad\times\HG{{}_2 F_1}{1-q_0,\,\, \frac{3}{4}-q_0}{2 - 2q_0}{-8 v (1 + v) (1 + 2 v + 2 v^2)}
\nonumber\\
&\quad\times\HG{{}_2 F_1}{q_0 - \frac{1}{4},\,\, q_0}{2q_0}
{8 u v (1 - u v) (1 - 2 u v + 2 u^2 v^2)},\label{eq:gaussproduct_right}
\end{align}
which is also holomorphic at $(u,v)=(0,0)$. Since the constant terms of \eqref{eq:H1leftside}
and \eqref{eq:gaussproduct_right} are equal these power series expansions must be equal.

\section{Horn's $H_5$}
\label{chap-H5}
The system of differential equations corresponding to Horn's $H_5$ function \eqref{eq-H5}, given in \cite[p.817]{Deb02}, is 
\begin{align}
&x(1+4x)F_{xx} -y(4x-1)F_{xy} - y^2F_{yy} \nonumber \\
&\,\,+(1-b+(4a-6)x)F_{x} + 2(a+1)yF_y + a(a+1)F = 0 \label{eq:H5_1} \\
&y(1-y)F_{yy} - xyF_{xy} + 2x^2F_{xx} \nonumber \\
&\,\,+(c-(a+b+1)y)F_y + (2+a-2b)xF_x -abF = 0.\label{eq:H5_2}
\end{align}


Again we would like to parameterize $x$ and $y$ with variable $t$ such that we get an equation of the form 
$$\frac{d^2 f}{d t^2} + c_1 \frac{d f}{d t} + c_2 f = 0.$$

First we can eliminate $F_{xx}$ and
$F_{xy}$ from \eqref{eq:Ftt}
using equations \eqref{eq:H5_1} and \eqref{eq:H5_2}. Then
we set the coefficient of 
$F_{yy}$ to be zero. This coefficient is equal to
{\small \begin{equation}\frac{(3 \, x y^{2} - 4 \, x y - y^{2} + y)\dot{x}^2 + (-12 \, x^{2} y + 8 \, x^{2} - 2 \, x y + 2 \, x)\dot{x}\dot{y} + (12 \, x^{3} - x^{2})\dot{y}^2}{{\left(12 \, x - 1\right)} x^{2}}.\label{eq:H5_coef}\end{equation}}

Instead of factoring \eqref{eq:H5_coef} such that we only have to solve linear differential equations. It is easier to solve this right away, thanks to a suggestion from Wadim Zudilin. Namely let $x$ be dependent on $y$ and try to find an algebraic relation between $x$ and $y$ using power series solutions of $x(y)$. Rewrite the numerator of \eqref{eq:H5_coef}, then we want to solve

\begin{equation}(3 \, x y^{2} - 4 \, x y - y^{2} + y)\frac{d^2x}{dy^2} + (-12 \, x^{2} y + 8 \, x^{2} - 2 \, x y + 2 \, x)\frac{dx}{dy} + (12 \, x^{3} - x^{2}) = 0.\label{eq:H5_diff}\end{equation}
Now by trying out different constant terms for $x(y)$ we can generate power series solutions from the recurrences of \eqref{eq:H5_diff} upto a certain degree $M$. Suppose we have found power series solutions $x_1(y),\ldots,x_N(y)$ all of degree $M$. Then we can make a matrix $Q_i$ where each row corresponds to the coefficient vector of $y^vx_i^u + O(y^M)$ and where the rows run over a sufficient number of pairs $(u,v) \in \mathbb{Z}_{\geq 0}^2$.

In our case it was enough to let $M=30$, the constant terms of $x(y)$ were chosen to be the twenty integers $\{-10,\ldots,-1\}\cup\{1,\ldots,10\}$, and we let $(u,v) \in \{(a,b): a = 0,\ldots,3,\,\,b=0,\ldots,2\}$. This is because we don't expect the algebraic relation between $x$ and $y$ to be very complex.

The left kernel of each $Q_i$ now corresponds to an algebraic relations between $x$ and $y$. By our choice of $x_1(y),\ldots,x_{20}(y)$ and the choice of $(u,v)$ this kernel turns out to be $1$-dimensional for each of these $Q_i$. Let the algebraic relation that generates the left kernel of $Q_i$ be denoted by $f_i$. 

Then we interpolate the coefficients of $\{f_i\}_{i=1,\ldots,N}$ separately using variable $a$ to form an algebraic relation \begin{equation}f(a,x,y) := \sum_{u=0}^3\sum_{v=0}^2 c_{u,v}(a)x^uy^v.\label{eq-alg-h5}\end{equation}
The full relation is given in Appendix \ref{sec-app-h5}.
Solving $f(a,x,y) = 0$ for $y$ gives two solutions

{\footnotesize
$$
y_1 = \frac{4 \, {\left(144 \, a^{2} x - 144 \, a^{\frac{3}{2}} x^{\frac{3}{2}} + 4 \, a^{2} - 36 \, a^{\frac{3}{2}} \sqrt{x} + 36 \, a x - 4 \, \sqrt{a} x^{\frac{3}{2}} + a - \sqrt{a} \sqrt{x}\right)} {\left(4 \, a + 1\right)}}{{\left(12 \, a - 1\right)}^{3} x},
$$
$$y_2 = \frac{4 \, {\left(144 \, a^{2} x + 144 \, a^{\frac{3}{2}} x^{\frac{3}{2}} + 4 \, a^{2} + 36 \, a^{\frac{3}{2}} \sqrt{x} + 36 \, a x + 4 \, \sqrt{a} x^{\frac{3}{2}} + a + \sqrt{a} \sqrt{x}\right)} {\left(4 \, a + 1\right)}}{{\left(12 \, a - 1\right)}^{3} x}.$$
}

Pick $y=y_2$ and substitute $x \to t^2$ and $a \to s^2$ to get the specialization
$$x(s,t) = t^2,\,\,y(s,t) = \frac{4 \, {\left(144 s^{2} t^{2} + 4 s^{2} + 32 s t + 4 t^{2} + 1\right)} {\left(4 \, s^{2} + 1\right)} {\left(s + t\right)} s}{{\left(12 \, s^{2} - 1\right)}^{3} t^{2}}.$$
One checks that \eqref{eq:H5_coef} is annihilated by these choices of $x$ and $y$.

This means in \eqref{eq:Ftt} we eliminated the coefficients $F_{xx}, F_{xy}$ and $F_{yy}$, and we would like to eliminate $F_x$ and $F_y$ using \eqref{eq:Ft}. Again we call the coefficients for $F_x$ and $F_y$ respectively $c_3$ and $c_4$ and the constant coefficient $c_5$. These coefficients can be found in Appendix \ref{sec-app-h5}. And miraculously we find $c_4 / c_3 = \dot{y} /\dot{x}$ when $a = q, b = q - \frac12, c = 2q$. This means that for fixed $s$ the function
$$\HG{H_5}{q,\,\, q-\frac12}{2q}{t^2,\frac{4 \, {\left(144 \, s^{2} t^{2} + 4 \, s^{2} + 32 \, s t + 4 \, t^{2} + 1\right)} {\left(4 \, s^{2} + 1\right)} {\left(s + t\right)} s}{{\left(12 \, s^{2} - 1\right)}^{3} t^{2}}}.$$
satisfies an ordinary Fuchsian equation of the form
$$\frac{d^2F}{dt^2} + r_2 \frac{dF}{dt} + r_3 F.$$
The coefficients $r_2$ and $r_3$ can be found in Appendix \ref{sec-app-h5}. The local exponents for this differential equation look like
\begin{center}
\begin{tabular}{c|cc}
Singularity & Exponent 1 & Exponent 2 \\ \hline
$t = 0$ & $2q - 1$ & $2q$ \\
$t = \infty$ & $q$ & $q+1$ \\
$t = -s$ & $0$ & $1-2q$ \\
Roots of $144s^2t^2 + 4s^2 + 32st + 4t^2 + 1 = 0$ & $0$ & $1-2q$
\end{tabular}
\end{center}
Consider the covering $$z \to -\frac{4(s+t)^2}{(12st+1)^2}.$$
The singularity $t=-s$ now corresponds to $z=0$ with multiplicity $2$. The singularity at a root of $144s^2t^2 + 4s^2 + 32st + 4t^2 + 1 = 0$ now corresponds to $z=1$ with multiplicity $1$. And lastly $z = \infty$ corresponds to the regular point $t = -\frac{1}{12s}$ with multiplicity $2$. This means the local exponents of ${}_2F_1(z)$ under this covering will look like
\begin{center}
\begin{tabular}{c|cc}
Singularity & Exponent 1 & Exponent 2 \\ \hline
$t = 0$ & $0$ & $1$ \\
$t = \infty$ & $0$ & $1$ \\
$t = -s$ & $0$ & $2(1-c)$ \\
Roots of $144s^2t^2 + 4s^2 + 32st + 4t^2 + 1 = 0$ & $0$ & $c-a-b$ \\
$t = -\frac{1}{12s}$ & 2a & 2b
\end{tabular}
\end{center}
Comparing the two tables we want to solve the equations $2(1-c) = 1-2q, 2(b-a)=1$ and $c-a-b=1-2q$. Which comes down to
$$a = \frac{3}{2}q - \frac{1}{2},\,\,b = \frac{3}{2}q,\,\,q+\frac12.$$
Making the last table
\begin{center}
\begin{tabular}{c|cc}
Singularity & Exponent 1 & Exponent 2 \\ \hline
$t = 0$ & $0$ & $1$ \\
$t = \infty$ & $0$ & $1$ \\
$t = -s$ & $0$ & $1-2q$ \\
Roots of $144s^2t^2 + 4s^2 + 32st + 4t^2 + 1 = 0$ & $0$ & $1-2q$ \\
$t = -\frac{1}{12s}$ & 3q-1 & 3q
\end{tabular}
\end{center}
Hence we need to multiply by a factor $(12st + 1)^{1-3q}t^{2q-1}$ to make the tables for $H_5$ and ${}_2F_1$ equal.
As a consequence we get that
\begin{equation}\HG{H_5}{q,\,\, q-\frac12}{2q}{t^2,\frac{4 \, {\left(144 \, s^{2} t^{2} + 4 \, s^{2} + 32 \, s t + 4 \, t^{2} + 1\right)} {\left(4 \, s^{2} + 1\right)} {\left(s + t\right)} s}{{\left(12 \, s^{2} - 1\right)}^{3} t^{2}}}.\label{eq:H5_mid}\end{equation}
satisfies the same differential equation as the one coming from
\begin{equation}(12st + 1)^{1-3q}t^{2q-1}\HG{{}_2F_1}{\frac32q-\frac12,\,\,\frac32q}{q+\frac12}{-\frac{4(s+t)^2}{(12st+1)^2}}.\label{eq:H5_mid_2F1}\end{equation}

Next we want the arguments of $H_5$ in \eqref{eq:H5_mid} to be symmetric, i.e. we want to find a parameterization of $s(u,v)$ and $t(u,v)$ such that
{\small
$$
H_5(\phi(s(u,v),t(u,v)),\psi(s(u,v),t(u,v))) = H_5(\phi(s(v,u),t(v,u)),\psi(s(v,u),t(v,u))).
$$
}
And we want the argument of ${}_2F_1$ in \eqref{eq:H5_mid_2F1} to be independent of $s$. If we can achieve this then we can swap roles of $u$ and $v$ to make the Bailey identity by a symmetry argument. There are two approaches to this. The first approach is to make $H_5$ symmetric by plugging in a power series for $t$ with variable coefficients. Then determine these coefficients and try to figure out which rational function belongs to this power series. When we perform this calculation where we fix $s(u,v) = \frac{u}{2}$, we find that by a miracle the argument of ${}_2F_1$ becomes independent of $s$. This approach is succesful but the result suggests a more intuitive approach which I will present here.

We want to find the inverse function of
$$g_s(t) = \frac{2(s+t)}{12st+1}.$$
This is simply
$$g_s^{-1}(t) = \frac{2t - s}{12st - 2}.$$
Hence if we take $s(u,v) = \frac{u}{2}$ and $t(u,v) = g_{u/2}^{-1}(v)$ this means \eqref{eq:H5_mid_2F1} becomes

\begin{equation}\left(3uv+1\right)^{q}
\left(1 - 3 \, u^{2}\right)^{1-3q}
\left(-\frac{u + v}{2}\right)^{2q-1}
{}_2F_1\left(\begin{array}{c|}\frac{3q-1}{2},\,\, \frac{3q}{2} \\ q + \frac{1}{2} \end{array}\,\,-v^2\right).\label{eq:H5_presym}\end{equation}

Miraculously \eqref{eq:H5_mid} becomes symmetric under this parameterization,
\begin{equation}H_5\left(\begin{array}{c|}
q,\,\, q - \frac{1}{2} \\ 2q
\end{array}\,\,
\frac{{\left(u + v\right)}^{2}}{4 \, {\left(3 \, u v + 1\right)}^{2}},
\frac{4 \, {\left(u^{2} + 1\right)} {\left(v^{2} + 1\right)} u v}{{\left(3 \, u v + 1\right)} {\left(u + v\right)}^{2}}
\right).\label{eq:H5_sym}\end{equation}

Multiply \eqref{eq:H5_presym} by the constant $$(-2)^{1-2q}(1-3u^2)^{3q-1}{}_2F_1\left(\begin{array}{c|}\frac{3q-1}{2},\,\, \frac{3q}{2} \\ q + \frac{1}{2} \end{array}\,\,-u^2\right).$$
This means that \eqref{eq:H5_presym} turns into something symmetric: 
$$\left(3uv+1\right)^{q}
\left(-\frac{u + v}{2}\right)^{2q-1}
{}_2F_1\left(\begin{array}{c|}\frac{3q-1}{2},\,\, \frac{3q}{2} \\ q + \frac{1}{2} \end{array}\,\,-u^2\right){}_2F_1\left(\begin{array}{c|}\frac{3q-1}{2},\,\, \frac{3q}{2} \\ q + \frac{1}{2} \end{array}\,\,-v^2\right).$$
By this symmetry we note that it satisfies the same differential equations as \eqref{eq:H5_sym} in both $u$ and $v$.

Consider the function
$$\eta: (u,v) \mapsto \left(\frac{{\left(u + v\right)}^{2}}{4 \, {\left(3 \, u v + 1\right)}^{2}},
\frac{4 \, {\left(u^{2} + 1\right)} {\left(v^{2} + 1\right)} u v}{{\left(3 \, u v + 1\right)} {\left(u + v\right)}^{2}}\right).$$
We want to pick a point $(u_0,v_0)$ with $-u_0^2,-v_0^2 \in \{0,1\}$, such that $\eta$ maps an open neighborhood of $(u_0,v_0)$ to an open neighborhood of $(0,0)$. Just as we saw with $H_1$, this is impossible. To get a meaningful Bailey identity we will restrict $\eta$ to a neighborhood of $(0,0)$ of the form $u = x, v = -x\cdot y$. By comparing the solution spaces we find and making both constant terms equal we see that
$$
\HG{H_5}{q,\,\,q-\frac{1}{2}}{2q}{\frac{x^2(y-1)^2}{4(3x^2y-1)^2},\frac{4(x^2y^2 + 1)(x^2+1)y}{(3x^2y-1)(y-1)^2}}=
$$
$$
(1-3x^2y)^q(1-y)^{2q-1}\HG{{}_2F_1}{\frac32q - \frac12,\,\,\frac32q}{q+\frac12}{-x^2y^2}\HG{{}_2F_1}{\frac12q,\,\,\frac12q+\frac12}{\frac32-q}{-x^2}.
$$
Note that this identity only depends on $x^2$ and $y$. So we may substitute $x \to \sqrt{x}$ to obtain the Bailey-type identity

\begin{equation}
\HG{H_5}{q,\,\,q-\frac{1}{2}}{2q}{\frac{x(y-1)^2}{4(3xy-1)^2},\frac{4(xy^2 + 1)(x+1)y}{(3xy-1)(y-1)^2}}=
$$
$$
(1-3xy)^q(1-y)^{2q-1}\HG{{}_2F_1}{\frac32q - \frac12,\,\,\frac32q}{q+\frac12}{-xy^2}\HG{{}_2F_1}{\frac12q,\,\,\frac12q+\frac12}{\frac32-q}{-x}.
\end{equation}
\newpage
\begin{appendices}
\section{Additional Coefficients}
\label{app-coef}
In sections \ref{chap-H4},\ref{chap-H1} and \ref{chap-H5} we left out some coefficients to improve readability. These will be displayed here.
\subsection{Horn's $H_4$}
\label{sec-app-h4}
The following coefficients belong to those of equation \eqref{eq-c3c4c5-form}. 
\begin{align*}
    c_3 &= \left(2 \, b s^{4} t^{6} + c s^{4} t^{6} + 2 \, a s^{3} t^{7} + b s^{3} t^{7} + a s^{2} t^{8} - 4 \, b s^{4} t^{5} - 6 \, b s^{3} t^{6} - 2 \, b s^{2} t^{7} \right.\\ &\left. + 2 \, s^{3} t^{7} + s^{2} t^{8} + 4 \, b s^{4} t^{4} - 2 \, c s^{4} t^{4} - 2 \, a s^{3} t^{5} + 15 \, b s^{3} t^{5} + 4 \, a s^{2} t^{6} + 8 \, b s^{2} t^{6}\right.\\ &\left. - 2 \, c s^{2} t^{6} + 2 \, a s t^{7} + b s t^{7} - 4 \, b s^{4} t^{3} - 20 \, b s^{3} t^{4} - 2 \, s^{4} t^{4} - 22 \, b s^{2} t^{5} - 5 \, s^{3} t^{5} \right.\\ &\left. - 6 \, b s t^{6} + 3 \, s^{2} t^{6} + 2 \, s t^{7} + 2 \, b s^{4} t^{2} + c s^{4} t^{2} - 2 \, a s^{3} t^{3} + 15 \, b s^{3} t^{3} - 10 \, a s^{2} t^{4} \right.\\ &\left. + 32 \, b s^{2} t^{4} + 4 \, c s^{2} t^{4} - 2 \, a s t^{5} + 15 \, b s t^{5} + 2 \, b t^{6} + c t^{6} - 6 \, b s^{3} t^{2} - 2 \, s^{4} t^{2} \right.\\ &\left. - 22 \, b s^{2} t^{3} - 12 \, s^{3} t^{3} - 20 \, b s t^{4} - 21 \, s^{2} t^{4}  - 4 \, b t^{5} - 5 \, s t^{5} + 2 \, a s^{3} t + b s^{3} t \right.\\ &\left. + 4 \, a s^{2} t^{2} + 8 \, b s^{2} t^{2} - 2 \, c s^{2} t^{2} - 2 \, a s t^{3} + 15 \, b s t^{3} + 4 \, b t^{4} - 2 \, c t^{4} - 2 \, b s^{2} t \right.\\ &\left. - s^{3} t - 6 \, b s t^{2} - 7 \, s^{2} t^{2} - 4 \, b t^{3} - 12 \, s t^{3} - 2 \, t^{4} + a s^{2} + 2 \, a s t + b s t + 2 \, b t^{2} \right.\\ &\left. + c t^{2} - s t - 2 \, t^{2}\right) s/\left({\left(s t^{2} + s + 2 \, t\right)} {\left(2 \, s t + t^{2} + 1\right)} {\left(s^2 - 1\right)}^{2} t^{4}\right)\\
\end{align*}
\begin{align*}
    c_4 &= -\left(a s^{2} t^{5} + b s^{2} t^{5} + a s t^{6} - 2 \, d s^{2} t^{4} - 2 \, b s t^{5} + s^{2} t^{5} + s t^{6} - 2 \, a s^{2} t^{3} \right.\\ &\left. - 2 \, b s^{2} t^{3} + 4 \, d s^{2} t^{3} - a s t^{4} + 4 \, d s t^{4} + a t^{5} + b t^{5} - 2 \, d s^{2} t^{2} + 4 \, b s t^{3} - 8 \, d s t^{3} \right.\\ &\left. - 4 \, s^{2} t^{3} - 2 \, d t^{4} - 2 \, s t^{4} + t^{5} + a s^{2} t + b s^{2} t - a s t^{2} + 4 \, d s t^{2} - 2 \, a t^{3} - 2 \, b t^{3} \right.\\ &\left. + 4 \, d t^{3} - 2 \, b s t - s^{2} t - 2 \, d t^{2} - 7 \, s t^{2} - 4 \, t^{3} + a s + a t + b t - t\right)\\
    &/\left({\left(s t^{2} + s + 2 \, t\right)} {\left(2 \, s t + t^{2} + 1\right)} t^{3}\right)\\
    c_5 &= -\frac{{\left(2 \, b s^{2} t + a s t^{2} + 2 \, a s t - 4 \, b s t + s t^{2} + a s + 2 \, b t + 2 \, s t + s\right)} a {\left(t - 1\right)}^{2}}{{\left(s t^{2} + s + 2 \, t\right)} {\left(2 \, s t + t^{2} + 1\right)} t^{2}}.\\
\end{align*}
The following coefficients belong to those of equation \eqref{eq:Fuchs2}.
\begin{align*}
r_2 &= 2 \, \left(q_{0} s^{2} t^{5} + q_{1} s^{2} t^{5} + q_{0} s t^{6} - 4 \, q_{1} s^{2} t^{4} - 2 \, q_{1} s t^{5} + s^{2} t^{5} + s t^{6} - 2 \, q_{0} s^{2} t^{3}  \right.\\ &\left. + 6 \, q_{1} s^{2} t^{3} - q_{0} s t^{4} + 8 \, q_{1} s t^{4} + q_{0} t^{5} + q_{1} t^{5} - 4 \, q_{1} s^{2} t^{2} - 12 \, q_{1} s t^{3} - 4 \, s^{2} t^{3} \right.\\ &\left. - 4 \, q_{1} t^{4} - 2 \, s t^{4} + t^{5} + q_{0} s^{2} t + q_{1} s^{2} t - q_{0} s t^{2} + 8 \, q_{1} s t^{2} - 2 \, q_{0} t^{3} + 6 \, q_{1} t^{3}  \right.\\ &\left. - 2 \, q_{1} s t - s^{2} t - 4 \, q_{1} t^{2} - 7 \, s t^{2} - 4 \, t^{3} + q_{0} s + q_{0} t + q_{1} t - t\right)\\
&/\left({\left(s t^{2} + s + 2 \, t\right)} {\left(2 \, s t + t^{2} + 1\right)} {\left(t + 1\right)} {\left(t - 1\right)} t\right) \\
r_3 &=  \frac{{\left(2 \, q_{1} s^{2} t + q_{0} s t^{2} + 2 \, q_{0} s t - 4 \, q_{1} s t + s t^{2} + q_{0} s + 2 \, q_{1} t + 2 \, s t + s\right)} q_{0} {\left(t - 1\right)}^{2}}{{\left(s t^{2} + s + 2 \, t\right)} {\left(2 \, s t + t^{2} + 1\right)} t^{2}}.
\end{align*}

\subsection{Horn's $H_1$}
\label{sec-app-h1}
Similarly we show the coefficients corresponding to the same equations, but now for Horn's $H_1$.
\begin{align*}
c_3 &= 2 \, \left(8 \, c s^{4} t^{6} + 2 \, d s^{4} t^{6} + 4 \, a s^{3} t^{7} + 4 \, b s^{3} t^{7} + 12 \, c s^{3} t^{7} + 2 \, a s^{2} t^{8} + 2 \, b s^{2} t^{8}\right.\\ &\left. + 4 \, c s^{2} t^{8} + 2 \, s^{3} t^{7} + s^{2} t^{8} + 8 \, c s^{4} t^{4} - 4 \, d s^{4} t^{4} - 4 \, a s^{3} t^{5} - 4 \, b s^{3} t^{5} + 40 \, c s^{3} t^{5}\right.\\ &\left. + 8 \, a s^{2} t^{6} + 8 \, b s^{2} t^{6} + 44 \, c s^{2} t^{6} - 4 \, d s^{2} t^{6} + 4 \, a s t^{7} + 4 \, b s t^{7} + 12 \, c s t^{7} - 4 \, s^{4} t^{4}\right.\\ &\left. - 8 \, s^{3} t^{5} + 2 \, s^{2} t^{6} + 2 \, s t^{7} + 2 \, d s^{4} t^{2} - 4 \, a s^{3} t^{3} - 4 \, b s^{3} t^{3} + 12 \, c s^{3} t^{3} - 20 \, a s^{2} t^{4}\right.\\ &\left. - 20 \, b s^{2} t^{4} + 44 \, c s^{2} t^{4} + 8 \, d s^{2} t^{4} - 4 \, a s t^{5} - 4 \, b s t^{5} + 40 \, c s t^{5} + 8 \, c t^{6} + 2 \, d t^{6}\right.\\ &\left. - 4 \, s^{4} t^{2}- 22 \, s^{3} t^{3} - 32 \, s^{2} t^{4} - 8 \, s t^{5} + 4 \, a s^{3} t + 4 \, b s^{3} t + 8 \, a s^{2} t^{2} + 8 \, b s^{2} t^{2}\right.\\ &\left.  + 4 \, c s^{2} t^{2} - 4 \, d s^{2} t^{2} - 4 \, a s t^{3} - 4 \, b s t^{3} + 12 \, c s t^{3} + 8 \, c t^{4} - 4 \, d t^{4} - 4 \, s^{3} t - 18 \, s^{2} t^{2} \right.\\ &\left. - 22 \, s t^{3} - 4 \, t^{4} + 2 \, a s^{2} + 2 \, b s^{2} + 4 \, a s t + 4 \, b s t + 2 \, d t^{2} - s^{2} - 4 \, s t - 4 \, t^{2}\right) s \\
&/\left({\left(s t^{2} + s + 2 \, t\right)} {\left(2 \, s t + t^{2} + 1\right)} {\left(s + 1\right)}^{2} {\left(s - 1\right)}^{2} t^{4}\right)\\
c_4 &= -2 \, \left(4 \, b s^{2} t^{3} + 4 \, c s^{2} t^{3} + 2 \, a s t^{4} + 2 \, b s t^{4} + 4 \, c s t^{4} + 2 \, s^{2} t^{3} + s t^{4} - 4 \, a s^{2} t \right.\\ &\left.  - 8 \, a s t^{2} + 8 \, b s t^{2} + 4 \, c s t^{2} + 4 \, b t^{3} + 4 \, c t^{3} + 2 \, s^{2} t + 6 \, s t^{2} + 2 \, t^{3} \right.\\ &\left. - 2 \, a s - 2 \, b s - 4 \, a t + s + 2 \, t\right)/\left({\left(s t^{2} + s + 2 \, t\right)} {\left(2 \, s t + t^{2} + 1\right)}\right)\\
c_5 &= -\frac{4 \, {\left(2 \, c s^{2} t^{3} + a s t^{4} + 2 \, c s t^{4} - 2 \, a s t^{2} + 2 \, c s t^{2} + 2 \, c t^{3} + a s\right)} b}{{\left(s t^{2} + s + 2 \, t\right)} {\left(2 \, s t + t^{2} + 1\right)} t^{2}}
\end{align*}
\begin{align*}
r_2 &= 2 \, \left(2 \, q_{0} s^{2} t^{3} + 2 \, q_{0} s t^{4} + 2 \, s^{2} t^{3} + s t^{4} - 2 \, q_{0} s^{2} t + 2 \, q_{0} t^{3}\right.\\ &\left. + 2 \, s^{2} t + 6 \, s t^{2} + 2 \, t^{3} - 2 \, q_{0} s - 2 \, q_{0} t + s + 2 \, t\right)\\ &/\left(\left(s t^{2} + s + 2 \, t\right)\left(2 \, s t + t^{2} + 1\right) t\right) \\
r_3 &= \frac{2 \, {\left(2 \, q_{0} s t^{4} + 2 \, s^{2} t^{3} + s t^{4} - 4 \, q_{0} s t^{2} + 4 \, s t^{2} + 2 \, t^{3} + 2 \, q_{0} s - s\right)} q_{0}}{{\left(s t^{2} + s + 2 \, t\right)} {\left(2 \, s t + t^{2} + 1\right)} t^{2}} 
\end{align*}

\subsection{Horn's $H_5$}
\label{sec-app-h5}
The algebraic relation \eqref{eq-alg-h5} is given by
{\footnotesize
\begin{align*}f(a,x,y) = &{\left(2985984 a^{6} - 1492992 a^{5} + 311040 a^{4} - 34560  a^{3} + 2160 a^{2} - 72 a + 1\right)} x^{2} y^{2}\\
&- 256 \, {\left(20736 \, a^{5} + 11520 \, a^{4} + 1888 \, a^{3} + 80 \, a^{2} + a\right)} x^{3} \\
&- 288 \, {\left(27648 \, a^{6} + 6912 \, a^{5} - 1152 \, a^{4} - 160 \, a^{3} + 28 \, a^{2} - a\right)} x^{2} y \\
& + 128 \, {\left(41472 \, a^{6} + 20736 \, a^{5} + 4032 \, a^{4} + 704 \, a^{3} + 82 \, a^{2} - a\right)} x^{2}\\
&  - 8 \, {\left(27648 \, a^{6} + 6912 \, a^{5} - 1152 \, a^{4} - 160 \, a^{3} + 28 \, a^{2} - a\right)} x y \\
& + 16 \, {\left(18432 \, a^{6} - 2304 \, a^{5} - 4608 \, a^{4} - 736 \, a^{3} - 8 \, a^{2} - a\right)} x \\
&+ 4096 \, a^{6} + 4096 \, a^{5} + 1536 \, a^{4} + 256 \, a^{3} + 16 \, a^{2}.\end{align*}
}
\noindent The remaining coefficients appearing in section \ref{chap-H5} are given by

\begin{align*}
c_3 &= -2\left(192 \, a s^{3} t^{2} - 192 \, b s^{3} t^{2} + 288 \, a s^{2} t^{3} + 192 \, s^{3} t^{2} + 288 \, s^{2} t^{3} - 8 \, a s^{3} \right.\\&\left. - 8 \, b s^{3} - 72 \, b s^{2} t  + 48 \, a s t^{2} - 48 \, b s t^{2} + 8 \, a t^{3} + 4 \, s^{3} + 36 \, s^{2} t + 48 \, s t^{2} + 8 \, t^{3} \right.\\&\left. - 2 \, a s - 2 \, b s - 2 \, b t + s + t\right)/\left({\left(144 \, s^{2} t^{2} + 4 \, s^{2} + 32 \, s t + 4 \, t^{2} + 1\right)} {\left(s + t\right)}\right)
\end{align*}
\begin{align*}
c_4 &= 4\left( \left(13824 \, a s^{6} t^{4} + 27648 \, b s^{6} t^{4} - 20736 \, c s^{6} t^{4} - 27648 \, a s^{5} t^{5} + 27648 \, b s^{5} t^{5} \right.\right.\\&\left.\left. - 41472 \, a s^{4} t^{6} + 13824 \, s^{6} t^{4} - 27648 \, s^{5} t^{5} - 41472 \, s^{4} t^{6} - 1920 \, a s^{6} t^{2}  \right.\right.\\&\left.\left. - 1536 \, b s^{6} t^{2} + 1728 \, c s^{6} t^{2} + 4608 \, b s^{5} t^{3} - 3456 \, a s^{4} t^{4} + 13824 \, b s^{4} t^{4} \right.\right.\\&\left.\left. + 5184 \, c s^{4} t^{4} - 7680 \, a s^{3} t^{5} + 7680 \, b s^{3} t^{5} - 2304 \, a s^{2} t^{6} + 1536 \, s^{6} t^{2} \right.\right.\\&\left.\left. + 13824 \, s^{5} t^{3} + 6912 \, s^{4} t^{4} - 7680 \, s^{3} t^{5} - 2304 \, s^{2} t^{6} - 64 \, a s^{6} - 64 \, b s^{6} \right.\right.\\&\left.\left. - 576 \, a s^{5} t - 576 \, b s^{5} t - 960 \, a s^{4} t^{2} - 768 \, b s^{4} t^{2} - 432 \, c s^{4} t^{2} + 1280 \, b s^{3} t^{3} \right.\right.\\&\left.\left. + 288 \, a s^{2} t^{4} + 1728 \, b s^{2} t^{4} - 432 \, c s^{2} t^{4} - 192 \, a s t^{5} + 192 \, b s t^{5} - 32 \, a t^{6} \right.\right.\\&\left.\left. + 32 \, s^{6} + 576 \, s^{5} t + 3360 \, s^{4} t^{2} + 3840 \, s^{3} t^{3} + 864 \, s^{2} t^{4} - 192 \, s t^{5} - 32 \, t^{6} \right.\right.\\&\left.\left. - 32 \, a s^{4} - 32 \, b s^{4} - 160 \, a s^{3} t - 160 \, b s^{3} t - 120 \, a s^{2} t^{2} - 96 \, b s^{2} t^{2} + 36 \, c s^{2} t^{2} \right.\right.\\&\left.\left. + 32 \, b s t^{3} - 8 \, a t^{4} + 12 \, c t^{4} + 16 \, s^{4} + 160 \, s^{3} t + 240 \, s^{2} t^{2} + 96 \, s t^{3} \right.\right.\\&\left.\left. - 4 \, a s^{2} - 4 \, b s^{2} - 4 \, a s t - 4 \, b s t - c t^{2} + 2 \, s^{2} + 4 \, s t + 2 \, t^{2}\right)\left(4 \, s^{2} + 1\right) s\right)\\&/\left({{\left(144 \, s^{2} t^{2} + 4 \, s^{2} + 32 \, s t + 4 \, t^{2} + 1\right)} {\left(12 \, s^{2} - 1\right)}^{3} {\left(s + t\right)} t^{4}}\right)
\end{align*}

\begin{align*}
c_5 &= -4a\left(12 \, a s^{3} t^{2} - 48 b s^{3} t^{2} + 36 a s^{2} t^{3} + 12 s^{3} t^{2} + 36 s^{2} t^{3} \right.\\&\left. + 4 b s^{3} + 3 a s t^{2} - 12 b s t^{2} + a t^{3} + 3 s t^{2} + t^{3} + b s\right)\\&/\left(\left(144 \, s^{2} t^{2} + 4 \, s^{2} + 32 \, s t + 4 \, t^{2} + 1\right) \left(s + t\right) t^{2}\right)
\end{align*}

\begin{align*}
r_2 &= 2 \, \left(144 \, q s^{2} t^{3} + 144 \, s^{3} t^{2} + 144 \, s^{2} t^{3} - 8 \, q s^{3} - 36 \, q s^{2} t  \right.\\&\left. + 4 \, q t^{3} + 4 \, s^{3} + 36 \, s^{2} t + 36 \, s t^{2} + 4 \, t^{3} - 2 \, q s - q t + s + t\right)\\&/\left(\left(144 \, s^{2} t^{2} + 4 \, s^{2} + 32 \, s t + 4 \, t^{2} + 1\right) \left(s + t\right) t\right)
\end{align*}
\begin{align*}
r_3 &= -2 \, \left(\left(72 \, q s^{3} t^{2} - 72 \, q s^{2} t^{3} - 72 \, s^{3} t^{2} - 72 \, s^{2} t^{3} - 8 \, q s^{3} \right.\right.\\&\left.\left. + 18 \, q s t^{2} - 2 \, q t^{3} + 4 \, s^{3} - 18 \, s t^{2} - 2 \, t^{3} - 2 \, q s + s\right) q\right)\\&/\left(\left(144 \, s^{2} t^{2} + 4 \, s^{2} + 32 \, s t + 4 \, t^{2} + 1\right) \left(s + t\right) t^{2}\right).
\end{align*}

\end{appendices}
\bibliographystyle{alphaurl}
\bibliography{bib} 
\end{document}